\documentclass[11pt,oneside,reqno]{amsart}
\usepackage{amsmath,amssymb,amsthm,array,mathrsfs}
\usepackage[T1]{fontenc}
\usepackage[all]{xypic}
\usepackage{mathtools}
\usepackage[mathscr]{euscript}
\usepackage{graphicx}
\usepackage{extpfeil}

\usepackage[top=1in,bottom=1in,left=1in,right=1in]{geometry}
\usepackage{xcolor}
\usepackage{ tipa }
\usepackage{titlecaps}

\usepackage{hyperref}
\usepackage{cleveref}
\usepackage{lipsum}

\newcommand{\R}{\mathbb R}

\makeatletter
\newtheorem*{rep@theorem}{\rep@title}
\newcommand{\newreptheorem}[2]{%
\newenvironment{rep#1}[1]{%
 \def\rep@title{#2 \ref{##1}}%
 \begin{rep@theorem}}%
 {\end{rep@theorem}}}
\makeatother

\makeatletter
\newtheorem*{rep@corollary}{\rep@title}
\newcommand{\newrepcorollary}[2]{%
\newenvironment{rep#1}[1]{%
 \def\rep@title{#2 \ref{##1}}%
 \begin{rep@corollary}}%
 {\end{rep@corollary}}}
\makeatother

\newtheorem{theorem}{Theorem}[section]
\newtheorem*{theorem*}{Theorem*}

\newtheorem{lemma}[theorem]{Lemma}
\newtheorem{proposition}[theorem]{Proposition}
\newtheorem{remark}[theorem]{Remark}
\newtheorem{claim}[theorem]{Claim}

\newtheorem{question}[theorem]{Question}
\newreptheorem{theorem}{Theorem}
\newrepcorollary{corollary}{Corollary}

\crefname{theorem}{Theorem}{Theorems}
\crefname{lemma}{Lemma}{Lemmas}
\crefname{proposition}{Proposition}{Propositions}
\crefname{corollary}{Corollary}{Corollarys}
\crefname{section}{Section}{Sections}
\crefname{definition}{Definition}{Definitions}
\crefname{question}{question}{Questions}

\title{Towards Ivanov's meta-conjecture for Geodesic currents}
\author{Meenakshy Jyothis}
\date{}

\begin{document}

\maketitle

\begin{abstract} Given a closed, orientable surface $S$ of negative Euler characteristic, we study two automorphism groups: $Aut(\mathscr{C})$ and $Aut(\mathcal{ML})$, groups of homeomorphisms that preserve the intersection form in the space $\mathscr{C}$ of geodesic currents and the space $\mathcal{ML}$ of measured laminations. We prove that except in a few special cases, $Aut(\mathcal{ML})$ is isomorphic to the extended mapping class group. This theorem is a special case of \textit{Ivanov's meta-conjecture}. We investigate this question for $Aut(\mathscr{C})$. To demonstrate the difficulty in proving Ivanov's conjecture for $Aut(\mathscr{C})$, we construct infinite
family of pairs of closed curves that have the simple same marked length spectra and self intersection
number.
\end{abstract}

\section{Introduction}   

Let $S$ be a closed, orientable, finite type surface of genus $g \geq 2$ and let $Mod^{\displaystyle \pm}(S)$ denote the extended mapping class group of $S$. We consider $Aut(\mathcal{ML})$, the group of homeomorphisms of measured laminations on $S$ that preserve the intersection form. We prove that Ivanov's meta-conjecture holds for $Aut(\mathcal{ML})$. More precisely, we prove that except in a few special cases, $Aut(\mathcal{ML})$ is isomorphic to $Mod^{\displaystyle \pm}(S)$. 

Now for $Aut(\mathcal{C})$, the group of homeomorphisms of geodesic currents on $S$ that preserve the intersection form, we show that there is a surjective map $f:$ $Aut(\mathcal{C}) \rightarrow Mod^{\displaystyle \pm}(S)$. However it is not easy to see that this map is injective. To illustrate this we construct infinite family of closed curves $(\gamma_{n}, \gamma_{n}')$ that have the same simple marked length spectra and self intersection number. 

\subsection{Ivanov's meta conjecture for measured laminations}

Ivanov's meta conjecture states that when $g \geq 3$, then every object naturally associated to a surface $S$ that has a ``sufficiently rich'' structure has $Mod^{\displaystyle \pm}(S)$ as its group of automorphisms, and makes a similar claim for when $g = 2$ \cite{Iva06}. In his seminal work, Ivanov proved that the automorphism group of the curve complex of $S$ is $Mod^{\displaystyle \pm}(S)$ \cite{Iva97}. Ivanov's theorem inspired a number of results in the following years \cite{Irm06, IK07, Di12, AACLOSX, BM19}. Many of these results considered different complexes associated to a surface and showed that their automorphism group is $Mod^{\displaystyle \pm}(S)$, and their proofs use Ivanov's original theorem.

We consider the space of measured laminations on a surface, $\mathcal{ML}(S)$. Given a hyperbolic structure on $S$, a measured lamination is a closed subset of $S$ 
 foliated by geodesics, equipped with a transverse measure. Weighted simple closed multicurves are dense in $\mathcal{ML}(S)$. By work of Kerckhoff, the geometric intersection number on simple closed curves extends to an `intersection form' $i(\cdot , \cdot)$ on $\mathcal{ML}(S)$ \cite{Ker}. In particular, the intersection form is a continuous bilinear map $i(\cdot , \cdot): \mathcal{ML}(S) \times \mathcal{ML}(S) \rightarrow \mathbb{R}$, such that for any two simple closed curves $\gamma$ and $\delta$, the intersection form $i(\gamma, \delta)$ agrees with the their geometric intersection number.

In this paper we will be  considering the automorphism group $Aut(\mathcal{ML})$ that consists of homeomorphisms on $\mathcal{ML}(S)$ that preserve the intersection form. Let
\[
Aut(\mathcal{ML}) = \{ \phi: \mathcal{ML}(S) \rightarrow  \mathcal{ML}(S) \text{ homeo } | i(\lambda, \mu) = i(\phi(\lambda) , \phi(\mu)) \}
\]

We show that in most cases $Aut(\mathcal{ML})$ is isomorphic to the extended mapping class group.

\begin{theorem}\label{thm0.1}
Let $S_{g}$ be a closed, orientable, finite type surface of genus $g \geq 2$ and let $Aut(\mathcal{ML})$ denote the group of homeomorphisms on $\mathcal{ML}(S_{g})$ that preserve the intersection form. Then, for all $g \neq 2$
\[
Aut(\mathcal{ML}) \cong Mod^{\displaystyle \pm} (S_{g})
\]
For the surface of genus 2, 
\[
Aut(\mathcal{ML}) \cong \displaystyle{Mod^{ \pm} (S_{2})/H}
\]
Where $H$ is the order two subgroup generated by the hyperelliptic involution.   
\end{theorem}

\begin{remark}
 The above theorem is equivalent to the statement that $Aut(\mathcal{ML})$ is isomorphic to the automorphism group of the curve complex for all surface $S_{g}$ with $g \geq 2$.
 \end{remark} 

 After we finished writing, we discovered that Ken'ichi Ohshika and Athanase Papadopoulos also prove that $Aut(\mathcal{ML})$ is isomorphic to $Mod^{\pm}(S)$ \cite{OHSHIKA2018899}. However, they use quite different techniques from us. In particular, many of the results in this paper also apply to the general setting of geodesic currents. However, the Ivanov's meta conjecture for geodesic currents does not follow immediately from these results. Some obstructions we faced when generalizing the result to the setting of currents is discussed in the next section.

\subsection{The case of geodesic currents}\label{secigc} The space of geodesic currents on a surface, $\mathscr{C}(S)$, is an extension of $\mathcal{ML}$ that contains weighted \textit{non-simple} closed curves as a dense subset in the same way as $\mathcal{ML}$ contains weighted simple closed multicurves as a dense subspace. By work of Bonahon, the intersection form also extends to $\mathscr{C}(S)$ \cite{Bon86}. The space of geodesic currents comes equipped with the weak* topology.

Let $Aut(\mathscr{C})$ denote the group of homeomorphisms on $\mathscr{C}(S)$ that preserve the intersection form. It is already known that $\displaystyle{Mod^{\pm}(S)}$ embeds in $Aut(\mathscr{C})$. Because our maps preserve the intersection form, the homeomorphisms in $Aut(\mathscr{C})$ restrict to homeomorphisms in $Aut(\mathcal{ML})$.  For $g \geq 3$, \cref{thm0.1} gives us the following maps. 
\[
  \displaystyle{Mod^{\pm}(S) \hookrightarrow Aut(\mathscr{C}) \xtwoheadrightarrow{f} Aut(\mathcal{ML}) \xrightarrow{\cong} Mod^{\pm}(S)}  
\]

We are interested to know if Ivanov's theorem holds for $Aut(\mathscr{C})$. One reason why this is not immediate is because it is hard to show that the surjective map $f:Aut(\mathscr{C}) \twoheadrightarrow Aut(\mathcal{ML})$ is also injective. We construct infinite family of pairs of closed curves $(\gamma_{n}$, $\gamma_{n}')$ with the same self intersection number and simple marked length spectra defined in Section 2.2. The kernel of $f$ could contain an automorphism of currents that maps $\gamma_{n}$ to $\gamma_{n}'$.

\begin{theorem}
\label{thm0.3} Let $S$ be a surface of genus at least 2.
We can find infinitely many pairs of closed curves $\gamma_n$ and $\gamma_n'$ on $S$ such that $\gamma_n$ and $\gamma_n'$ have the same self intersection number and the same simple marked length spectra.
\end{theorem} 

There are results similar to \cref{thm0.3} in the literature. In fact, Parlier and Xu prove the existence of curves that satisfy an even stronger condition called $k$-equivalence. Two non-isotopic closed curves $\alpha$ and $\beta$ are said to be $k$-equivalent if they intersect all the closed curves with self intersection $k$ the same number of times \cite{kavi2019}. For any given $k > 0$, Parlier and Xu construct closed curves $\alpha$ and $\beta$ that are not $k$-equivalent but are $k'$-equivalent for any $k'< mk^2$ different from $k$ \cite{HpBx2023}.  In particular, the curves $\alpha$ and $\beta$ constructed are simple intersection equivalent. However, the curves $\alpha$ and $\beta$ have different self intersection number. This means that any map $\phi$ in the kernel of $f$ cannot map $\alpha$ to $\beta$. For the purpose of this paper we are interested in pairs of closed curves that share the same self intersection number as well.


\begin{remark}
    In the infinite family of closed curves that we construct, no two curves $(\gamma_{n}, \gamma_{n'})$ in the same pair are hyerbolically equivalent. In \cite{CL}, Leininger constructs examples of closed curves that are simple intersection equivalent but are not hyerbolically equivalent. Our curves add to these examples. In addition, the curves $\gamma_{n}$ and $\gamma_{n}'$ also have the same self intersection number.
\end{remark}

\subsection{Plan of the paper}
The paper is organized as follows. Section 2 provides background on measured laminations and geodesic currents. In section 3 we prove some properties of elements in $Aut(\mathcal{ML})$. The proof of \cref{thm0.1} follows immediately from these properties.
In particular, in section 3 we prove that any $\phi \in Aut(\mathcal{ML})$ is linear when $\mathcal{ML}$ is viewed as a subset of $\mathscr{C}(S)$. In this section, we also show that any such $\phi$ maps simple closed curves to weighted simple closed curves. As a consequence of these properties, we get that any $\phi \in Aut(\mathcal{ML})$ maps simple closed curves to simple closed curves, and therefore has an action on the curve complex. In section 4, we prove \cref{thm0.1} using Ivanov's original theorem and a density argument.

In section 5, we prove \cref{thm0.3} and we remark about some more properties of these curves. In this section we also discuss why these examples are a hindrance when extending Theorem 1.1 to the setting of geodesic currents.
In section 6, we discuss some open questions about $Aut(\mathcal{C})$.

\subsection{Acknowledgements.}
The author thanks her advisor Eugenia Sapir for her support throughout this project and for the many insightful conversations that greatly contributed to this work. The author also thanks Didac Martinez Granado for bringing Ken’ichi Ohshika's and Athanase Papadopoulos' work and Hugo Parlier's and Binbin Xu's work to their attention. 
\section{Background}

Let S be a hyperbolic surface with a complete metric defined by $\mathbb{H}^2/ \Gamma$, for $\Gamma \leq$ PSL(2,$\mathbb{R}$). We can identify the universal cover of the surface $S$ with $\mathbb{H}^2$. The spaces $\mathcal{ML}(S)$ and $\mathscr{C}(S)$, that were discussed in the introduction, are independent of the choice of hyperbolic metric on $S$ \cite{Mar2016}. However, in order to define a geodesic current or a measured lamination we will need to fix a hyperbolic metric.

Throughout the paper, subsurfaces of $S$ and closed curves will be considered up to isotopy. The complexity of a closed surface of genus $g$ is defined to be 3$g$ - 3.

\subsection{Laminations:} A \textit{geodesic lamination} $\lambda$ is a closed subset of $S$ foliated by simple, complete geodesics. An example of a geodesic lamination is a multicurve consisting of pairwise disjoint, simple closed geodesics on $S$. Another fundamental example arises from considering a set of disjoint geodesics in $\mathbb{H}^2$ that are invariant under $\Gamma$, whose union is closed. Projecting this set onto S gives us a geodesic lamination. We say a geodesic lamination is \textit{minimal} if it contains no proper non-empty sublamination.

 A \textit{complementary region} of a geodesic lamination $\lambda$ is a connected component of the open complement $S \backslash \lambda$. A geodesic lamination that \textit{fills} a subsurface $S'$ of $S$ intersects every essential, non-peripheral simple closed curve on $S'$. The complementary regions of such a filling lamination are either ideal polygons on $S'$ or crowns homotopic to a boundary curve of $S'$.

Some geodesic laminations can be equipped with a transverse measure, which assigns a positive measure to each arc $\tau$ that intersects the lamination $\lambda$ transversely. This measure is invariant under homotopy transverse to $\lambda$ and is supported on $\tau \cap \lambda$. A \textit{measured lamination} is a geodesic lamination with a transverse measure that has full support. Abusing notation, we will use $\lambda$ to denote the transverse measure. It is important to note that a measured lamination $\lambda$ can come equipped with different transverse measures. To avoid confusion between measures and their supports we have used notation more carefully later in this paper. A description of this can be found under `Notational Conventions' in \cref{NC}.

\subsection{Geodesic currents:}\label{sec2.2} Measured laminations can be thought of as a subset of a larger space of measures called geodesic currents defined as follows. Let $\mathcal{G}$ be the set of all unparameterized, unoriented complete geodesics in $\mathbb{H}^2$. Any geodesic in $\mathbb{H}^2$ can be determined by its extreme points on the boundary of $\mathbb{H}^2$. Hence there is a natural bijection  $\mathcal{G} \cong \displaystyle{( S^1 \times S^1 \setminus \Delta)/  \sim} $, where $\Delta$ is the diagonal in $S^1 \times S^1$ and the equivalence relation $\sim$ identifies the coordinates $(a,b)$ and $(b,a)$. We assign $\mathcal{G}$ the topology of $\displaystyle{( S^1 \times S^1 \setminus \Delta)/  \sim}$. The fundamental group, $\pi_1(S)$ embeds as a subgroup $\Gamma$ of PSL(2,$\R$). The isometry group PSL(2,$\R$) of $\mathbb{H}^2$ acts naturally on $\mathcal{G}$. A geodesic current is a $\Gamma$-invariant, locally finite, positive, Borel measure on $\mathcal{G}$. The space of all geodesic currents endowed with a weak* topology is denoted by $\mathscr{C}(S)$.

The set of closed curves on $S$ embeds into the space of currents. To see this consider a closed curve $\gamma$ on $S$. The lifts of $\gamma$ are a discrete subset of $\mathcal{G}$. The Dirac measure on this discrete set is a geodesic current,  and by abuse of notation we will denote the geodesic current by $\gamma$ as well.
Since positive measures are closed under addition and scalar multiplication by positive reals, weighted multicurves are also examples of geodesic currents. 

By work of Bonahon, it is known that weighted closed curves are dense in $\mathscr{C}(S)$ \cite{Bon86}. Bonahon also shows that geometric intersection number $i(.,.)$ defined on a pair of closed curve can be extended bilinearly and continuously to intersection form on a pair of currents.
For two any two currents $\mu$ and $\nu$, the intersection form $i(\mu, \nu)$ is defined as $\mu \times \nu (\mathcal{J}/ \Gamma)$, where $\mathcal{J}$ is a subset of $\mathcal{G} \times \mathcal{G}$ consisting of all pairs of incident distinct geodesics in $\mathbb{H}^2$.

Any geodesic current $\mu$ satisfying $i(\mu, \mu) = 0$ is supported on a $\Gamma$-invariant geodesic lamination on $\mathbb{H}^2$, and $\mu$ induces a transverse measure on its support. Because of this, $\mathcal{ML}$ embeds into $\mathscr{C}(S)$. 

For a geodesic current $\mu$, one can consider the intersection of $\mu$ with every closed curve on $S$. The infinite sequence $\big( i(\mu, \gamma_1), i(\mu, \gamma_2), \dots \big)$, where $\{\gamma\}_{i}$ is the collection of closed curves on $S$, is the marked length spectrum of $\mu$. Otal proved that every geodesic current $\mu$ can be identified using its unique marked length spectrum. \cite{Ot90}. The simple marked length spectrum of a current $\mu$ is the infinite sequence $\big( i(\mu, s_1), i(\mu, s_2),  \dots \big)$, where $\{s_i \}$ is the collection of simple closed curve on $S$.

\subsection{Notational Convention:}\label{NC} For a geodesic current $\mu$, consider the support of $\mu$ in $\mathcal{G}$ and take its projection onto $S$. We will denote this set by $|\mu|$. We will often use this notation to avoid confusion between the current and its support. For instance, the geodesic current corresponding to a closed curve $\gamma$ will be denoted by $\gamma$, but the curve itself will be denoted by $|\gamma|$. Similarly, the transverse measure on a measured lamination will be denoted by $\lambda$, and the lamination itself will be denoted by $|\lambda|$. For example, if a measured lamination is not uniquely ergodic and supports two distinct measures we will denote them using distinct symbols $\lambda$ and $\eta$. But, we will have $|\lambda| = |\eta|$.
\subsection{Curve Complex and Ivanov's theorem:}
The curve complex, $CC$ on a surface $S$ is a simplicial complex whose vertices correspond to isotopy classes of essential simple closed curves on $S$, and whose edges correspond to pairs of simple closed curves that have geometric intersection number zero. 

Let us denote the vertex set of the curve complex by $V(S)$. An automorphism of a curve complex is a bijection on $V(S)$ that takes simplices to simplices. We will use $Aut(CC)$ to denote the automorphism group of the curve complex on a surface $S$ . Ivanov's theorem states that $Aut(\mathcal{CC}) \cong Mod^{\displaystyle \pm} (S)$ for surfaces of genus at least three. For $S_2$, the closed, orientable surface of genus 2, $Aut(\mathcal{CC}) \cong \displaystyle{Mod^{ \pm} (S_{2})/H}$, where, $H$ denotes the subgroup generated by the hyperelliptic involution.

\section{Properties of automorphisms that preserve the intersection form} 

In this section we prove two properties of elements in $Aut(\mathcal{ML})$ that play a key role in proving $Aut(\mathcal{ML}) \cong Mod^{\displaystyle \pm} (S)$. Namely, we show that any $\phi \in Aut(\mathcal{ML})$ is linear and that any such $\phi$ maps simple closed curves to simple closed curves.

Most of the content in this section remains true in the setting of geodesic currents. In fact, the statements of Propositions 3.1, 3.3 and 3.13 hold true for any $\phi \in Aut(\mathscr{C})$. The remarks in this section concern about how the results generalizes to the space of currents. 

\subsection{Linearity}
\begin{proposition}\label{prop3.1}
Let $\phi \in Aut(\mathcal{ML})$ and let $\lambda , \nu \in \mathcal{ML}$ such that $\lambda + \nu \in \mathcal{ML}$. Let $c$ be a positive real number. Then
$$\phi( \lambda + \textup{c} \hspace{0.02cm} \nu) = \phi(\lambda) + \textup{c} \hspace{0.02cm} \phi(\nu) $$
\end{proposition}

\begin{proof} 
Observe that $Aut(\mathcal{ML})$ is a group, and if $\phi$ preserves the intersection form then $\phi^{-1}$ preserves the intersection form as well. Now, for any simple closed curve $\gamma$, we have the following equality
\begin{align*}
 i( \gamma, \phi(\lambda + \textup{c} \hspace{0.02cm} \nu)) & = i(\phi^{-1}(\gamma), \lambda + \textup{c} \hspace{0.02cm} \nu)\\
 & =  i(\phi^{-1}(\gamma), \lambda) + \textup{c} \hspace{0.02cm} i(\phi^{-1}(\gamma), \nu)\\
 & = i( \gamma, \phi(\lambda)) + \textup{c} \hspace{0.02cm} i(\gamma, \phi(\nu))\\
 & = i( \gamma, \phi(\lambda)) +  i(\gamma, \textup{c} \hspace{0.02cm} \phi(\nu))\\
 &  = i(\gamma, \phi(\lambda) + \textup{c} \hspace{0.02cm} \phi(\nu))
\end{align*}
This implies $\phi(\lambda + \textup{c} \hspace{0.02cm}\nu)$ and $\phi(\lambda) + \textup{c} \hspace{0.02cm}\phi(\nu)$ have the same simple marked length spectrum, and therefore $\phi(\lambda + \textup{c} \hspace{0.02cm}\nu) = \phi(\lambda) + \textup{c} \hspace{0.02cm}\phi(\nu)$ \cite{Primer}.
\end{proof}

\begin{remark} 
If $\gamma$ is allowed to be any closed curve, then the same proof can be used to see that any $\phi \in Aut(\mathscr{C})$ is linear. In this case, we will be looking at the marked length spectrum of a current instead of its simple marked length spectrum.  However, from Otal's work on currents it is known that a geodesic current is uniquely determined by its marked length spectrum~\cite{Ot90}.
\end{remark}

\subsection{Mapping simple closed curves to weighted simple closed curves}
\begin{proposition}\label{prop3.2}
Automorphisms of the space of measured laminations that preserve the intersection form map simple closed curves to weighted simple closed curves.
\end{proposition}

We prove this proposition after proving Lemmas 2.3 - 2.7.

\begin{lemma}\label{lma3.4}
Let $\lambda_1$ and $\lambda_2$ be any two geodesic currents satisfying $|\lambda_1| = |\lambda_2|$. Then for any geodesic current $\mu$, $i(\lambda_1, \mu ) \neq 0 \Longleftrightarrow i(\lambda_2, \mu) \neq 0$.

In particular, if $\lambda_1$ and $\lambda_2$ are two measured laminations satisfying $|\lambda_1| = |\lambda_2|$, then for any measured lamination $\mu$, $i(\lambda_1, \mu ) \neq 0 \Longleftrightarrow i(\lambda_2, \mu) \neq 0$.
\end{lemma}

\begin{proof} Let $\mathcal{J}$ be the set defined in \cref{sec2.2}. The support of $\lambda_{1} \times \mu$ in $\mathcal{J}$ is the set of ordered pair of geodesics $(g_{1}, g_{2})$ such that $g_{1} \in |\lambda_{1}|$ and $g_{2} \in |\mu|$ \cite{JACL}. That means, the support of $\lambda_{1} \times \mu$ in $\mathcal{J}$ is same as the support of $\lambda_{2} \times \mu$ in $\mathcal{J}$.  Therefore, the supports of $\lambda_{1} \times \mu$ and  $\lambda_{2} \times \mu$ in $\mathcal{J}/ \Gamma$ are the same. But then,
$i(\lambda_1, \mu ) \neq 0$ implies  $i(\lambda_2, \mu) \neq 0$ and vice versa. 
\end{proof}

The next lemma was previously known, see for example \cite{DS03}. But we include the proof here for completeness.
\begin{lemma} \label{lma3.5} 
Let $\lambda$ be a minimal measured lamination that fills a subsurface $S'$ of $S$ and let $\mu$ be any other measured lamination with non-empty support satisfying:
\begin{enumerate}
    \item $|\mu|$ is contained in the interior  of $S'$.
    \item $i(\lambda, \mu) = 0$
\end{enumerate}
Then $|\mu| = |\lambda|$.
\end{lemma}

\begin{proof}Let $\omega$ be a measured lamination consisting of leaves which are common to $|\lambda|$ and $|\mu|$. Since $\lambda$ is minimal, $|\omega|$ has to be either empty or is all of $|\lambda|$. 

\textbf{Case 1.}
If $\omega$ is empty, then the geodesics in $|\lambda|$ and $|\mu|$ neither coincide nor intersect. Also, observe that $|\mu|$ cannot intersect any of the  boundary curves of $S'$. This implies that, in the universal cover, any lift of $|\mu|$ has to live in the complementary regions formed by lifts of $|\lambda|$ and the lifts of the boundary curves of $S'$. As $\lambda$ is filling, such complementary regions will either be ideal polygons in $\mathbb{H}^2$ bounded by geodesics in the lift of $|\lambda|$ or they will be crowns bounded by both the geodesics in the lifts of $|\lambda|$ and the geodesics in the lifts of boundary curves of $S'$ (see Figure 1).

Any leaf in $\mu$ can only go from one ideal vertex of such an ideal polygons to another. But any such leaf is an isolated open leaf and cannot be in the support of a measured lamination \cite{CB88}. This contradicts our assumption that $\mu$ is a measured lamination with a non empty support. 

\begin{figure}
    \centering
    \includegraphics[width=\linewidth]{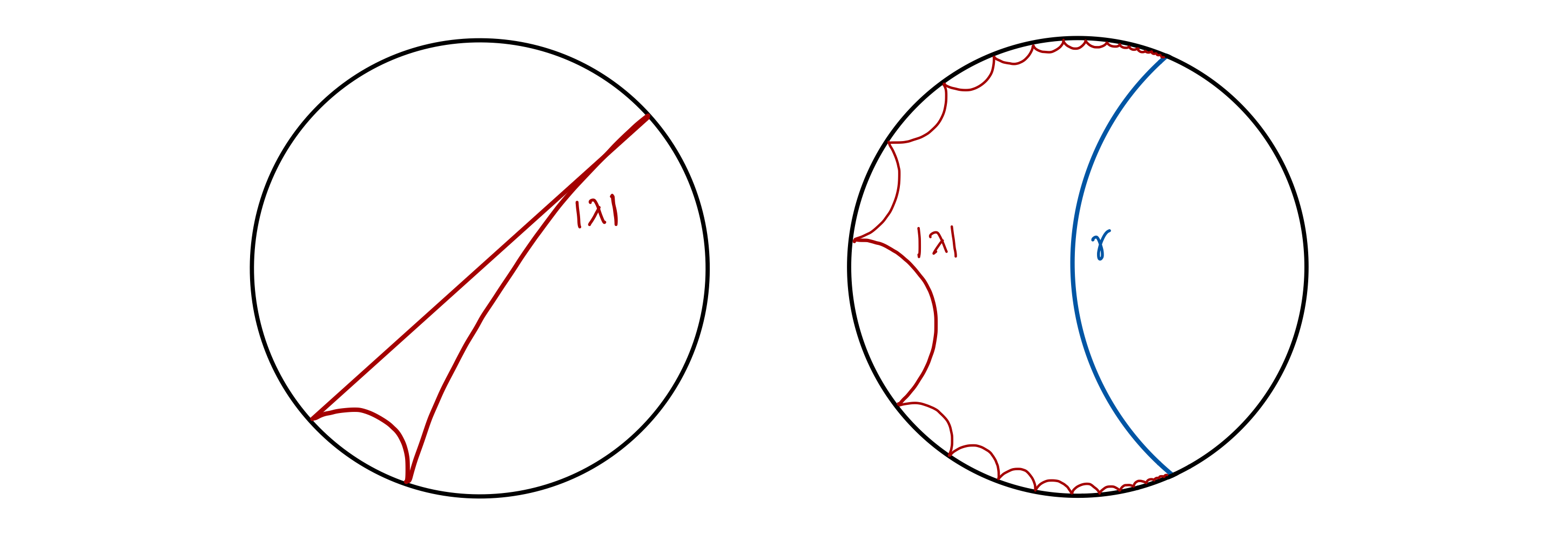}
    \caption{Two types of complimentary region of $|\lambda|$. Here, $\gamma$ is a peripheral curve of $S'$}
\end{figure}

\textbf{Case 2.} Now, if $\omega$ is all of $|\lambda|$, we get $|\lambda| = |\omega| \subseteq |\mu|$. Let $\mu_{\omega}$ denote the restriction of $\mu$ to $\omega$. Notably, $|\omega|$ is closed. Using decomposition of laminations into minimal components, we can write $\mu = \mu_{\omega} + \mu^{\prime}$. Here $\mu^{\prime}$ represents the measured lamination that encompasses the portion of $\mu$ disjoint from $\omega$. Observe that, both the conditions (1) and (2) in the statement of the lemma are fulfilled  by $\mu^{\prime}$. As we have already established in the proof of case 1, $\mu^{\prime}$ cannot support a measured lamination. Consequently, we conclude that $|\lambda|$ is equal to $|\mu|$.
\end{proof} 

Before proving more results, we want to define the set $\mathcal{E}_{\lambda}^{\mathcal{ML}}$, a generalized version of the set $\mathcal{E}_{\lambda}$ constructed in \cite{BIPP} \cite{BIPP2}. For a measured lamination $\lambda$, the set $\mathcal{E}_{\lambda}$ is defined as the set of all closed geodesics $c$ on $S$ that satisfy the following two conditions.

\begin{enumerate}
    \item $i(\lambda, c) = 0$
    \item For every closed curves $c'$ with $i(c, c') \neq 0$, 
     $i(\lambda, c') \neq 0$. 
\end{enumerate}

If $\lambda$ is a measured lamination that fills a subsurface $S'$ of $S$, then $\mathcal{E}_{\lambda}$ consists of all simple geodesics that forms the boundary of $S'$ \cite{BIPP}. We will call such curves peripheral curves of $S'$.

The set $\mathcal{E}_{\lambda}^{\mathcal{ML}}$ generalizes this construction to measured laminations and is defined as follows. For a measured lamination $\lambda$, the set $\mathcal{E}_{\lambda}^{\mathcal{ML}}$ is the set of all \textit{measured laminations} $\mu$ on $S$ that satisfy the following two conditions.

\begin{enumerate}
    \item $i(\lambda, \mu) = 0$
    \item For every measured lamination $\kappa$ with $i(\mu, \kappa) \neq 0$, $i(\lambda, \kappa) \neq 0$. 
\end{enumerate}

We make the following claims about $\mathcal{E}_{\lambda}^{\mathcal{ML}}$.
\begin{claim}
$\mathcal{E}_{\lambda}^{\mathcal{ML}}$ only consists of pairwise disjoint measured laminations. 
\end{claim}
\begin{proof}
   Any two measured laminations that intersect cannot simultaneously satisfy both conditions in the definition of $\mathcal{E}_{\lambda}^{\mathcal{ML}}$. 
   Let $\mu_{1}$ and $\mu_{2}$ be two measured laminations such that $\mu_{1} \in \mathcal{E}_{\lambda}^{\mathcal{ML}}$ and $i(\mu_1, \mu_2) \neq 0$. Then by the second condition in the definition of $\mathcal{E}_{\lambda}^{\mathcal{ML}}$, we get $i(\lambda, \mu_2) \neq 0$. Any such $\mu_2$ fails to belong to $\mathcal{E}_{\lambda}^{\mathcal{ML}}$.
\end{proof}
 
\begin{claim} The set $\mathcal{E}_{\lambda}^{\mathcal{ML}}$ is closed under addition and scalar multiplication by positive reals.
\end{claim}
\begin{proof}
    The proof readily follows from the bilinearity of the intersection form.
\end{proof}
Now, we will discuss what elements are contained in the set $\mathcal{E}_{\lambda}^{\mathcal{ML}}$ for some specific cases of $\lambda$.
\begin{lemma}\label{lma3.6}
If $\lambda$ is a simple closed curve, then $\mathcal{E}_{\lambda}^{\mathcal{ML}}$ consists of all the positive scalar multiples of $\lambda$. 
\[
\mathcal{E}_{\lambda}^{\mathcal{ML}} = \{c\lambda : c \in \mathbb{R}_{+}\}.
\]
\end{lemma}
\begin{proof} It is easy to see that all measured laminations $c \lambda$ with $c \in \mathbb{R}_{+}$ belong to $\mathcal{E}_{\lambda}^{\mathcal{ML}}$.

The set $\mathcal{E}_{\lambda}^{\mathcal{ML}}$ does not contain any other simple closed curve as they all fail to satisfy the second condition in the definition of $\mathcal{E}_{\lambda}^{\mathcal{ML}}$. For every simple closed curve $|\gamma| \neq |\lambda|$, we can find a simple closed curve that intersects $|\gamma|$ and not $|\lambda|$. 

The set $\mathcal{E}_{\lambda}^{\mathcal{ML}}$ also cannot contain any measured lamination that fills a subsurface $S'$ of $S$. Any subsurface $S'$  that $\mu$ fills will have complexity at least one. This means that we can find a curve $|\gamma|$ intersecting $\mu$ and contained entirely in $S'$. Now, if $\mu$ satisfies $i(\lambda, \mu) = 0$, then the curve $|\lambda|$ lies entirely outside of $S'$ and has zero intersection with $\gamma$. Therefore, $\mu$ fails to satisfy the second condition in the definition of $\mathcal{E}_{\lambda}^{\mathcal{ML}}$.
\end{proof}
\vspace{-0.5cm}
\begin{lemma}\label{lma3.7}
If $\lambda$ is a minimal measured lamination that fills a subsurface $S'$ of $S$, then the set $\mathcal{E}_{\lambda}^{\mathcal{ML}}$ is the smallest set that is closed under addition and scalar multiplication (by positive reals), that  consists of:
\begin{enumerate}
    \item peripheral curves that form the boundary of the subsurface $S'$,
    \item measured laminations supported on $|\lambda|$.
\end{enumerate}
\end{lemma}

\begin{proof}
The proof is divided into two parts. The first part shows that peripheral curves of $S'$ are the simple closed curves in $\mathcal{E}_{\lambda}^{\mathcal{ML}}$.
The second part focuses on measured laminations supported on $|\lambda|$. 

\textbf{Peripheral curves of} $\mathbf{S'}$\textbf{:}
To see the peripheral curves of $S'$ are the only simple closed curves in $\mathcal{E}_{\lambda}^{\mathcal{ML}}$, consider the complimentary subsurface $S \backslash S'$. We will denote this subsurface by $S''$. If $S''$ were a pair of pants, then it does not contain any other simple closed curve. And if it is a surface of a higher complexity, then for every simple closed curve $\delta$ in the interior of $S''$ we can find an another simple closed curve in the interior of $S''$ that also intersects $\delta$. So, simple closed  curve in the interior of $S''$ cannot satisfy both the conditions of $\mathcal{E}_{\lambda}^{\mathcal{ML}}$.  
Also, any closed curve contained entirely in $S'$ will intersect $\lambda$ and therefore cannot be contained in $\mathcal{E}_{\lambda}^{\mathcal{ML}}$.

Now, we show that all the peripheral curves of $S'$ belong to $\mathcal{E}_{\lambda}^{\mathcal{ML}}$. Consider $\gamma$, a peripheral curve of $S'$. As $\gamma \in \mathcal{E}_{\lambda}$ we know that $i(\lambda, \gamma) = 0$ and that if $\gamma$ intersects any closed curve then $\lambda$ intersects the same closed curve as well. It is left to check whether all the measured laminations that intersect $\gamma$ also intersect $\lambda$. To see that, consider the lift of $\gamma$ and $|\lambda|$ on $\mathbb{H}^2$. The complimentary region bounded by these geodesics will give us a crown as given in Figure 2.

\begin{figure}
    \centering
    \includegraphics[width=\linewidth]{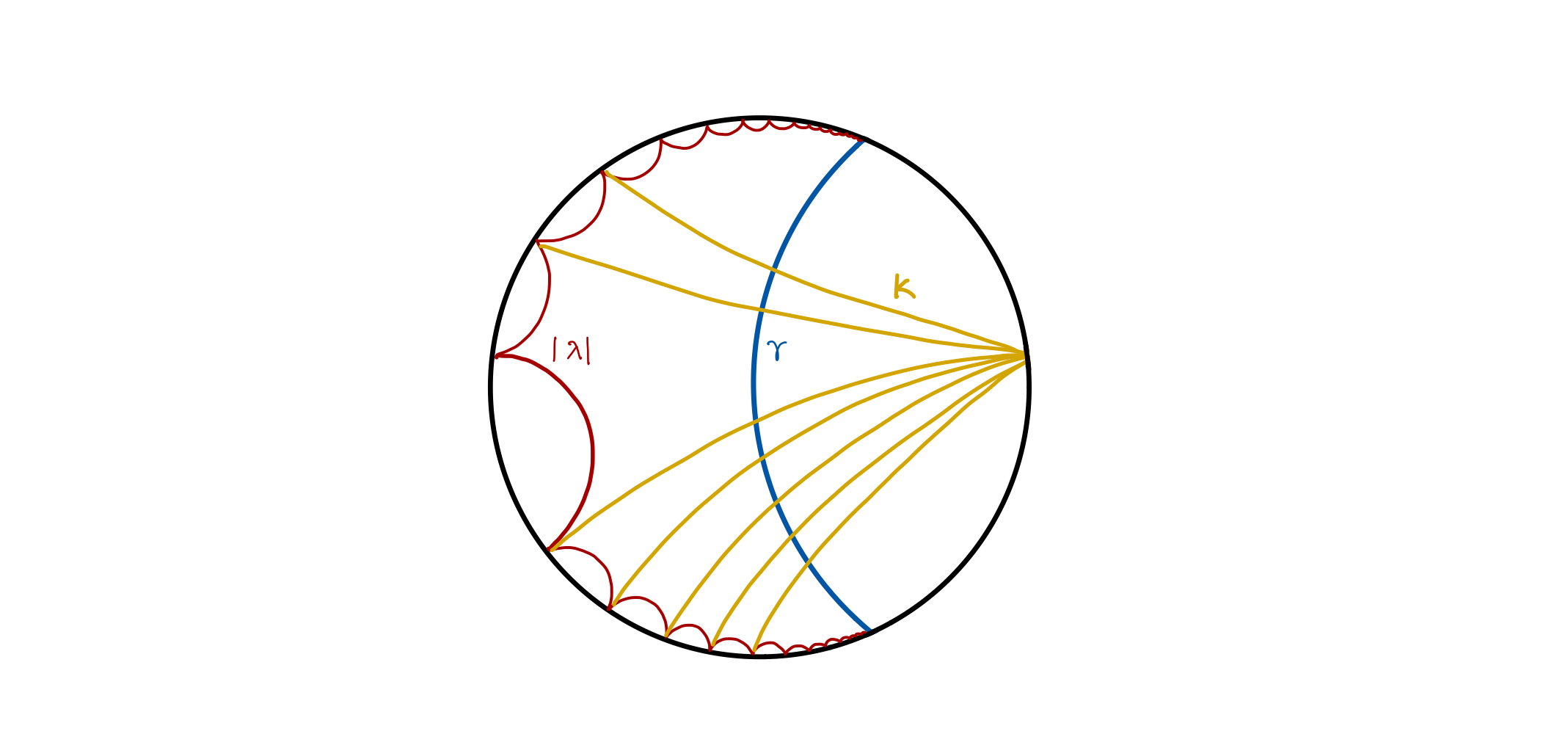}
    \caption{Example of a lamination $k$ that intersects $\gamma$ but not $\lambda$.}
\end{figure}

If there exist a lamination $k$ that intersects $\gamma$ but not $\lambda$, then the lift of $k$ will contain a leaf that has its one point in an ideal vertex of the crown, as given in Figure 2. But any such leaf will be an isolated leaf and cannot be in the support of a measured lamination  \cite{CB88}. Therefore, no such $k$ exist. 

 \textbf{Measured laminations supported on} $\mathbf{|\lambda|}$ \textbf{:} If $\mu$ is a measured lamination that fills a subsurface $S''$, and if $S'' \setminus S'$ is non empty, then $\mu$ will intersect a closed curve that does not intersect $\lambda$. This implies that any filling measured lamination in $\mathcal{E}_{\lambda}^{\mathcal{ML}}$ must have its support inside $S'$. But since these measured lamination also do not intersect $\lambda$, by~\cref{lma3.5} its support will be $|\lambda|$.
 \end{proof} 

  Since the criteria we used for defining $\mathcal{E}_{\lambda}^{\mathcal{ML}}$ only depends on the  intersection form, $\mathcal{E}_{\lambda}^{\mathcal{ML}}$ itself depends on only the support $|\lambda|$  . That is, if $|\lambda_{1}| = |\lambda_2| $, then $\mathcal{E}_{\lambda_1}^{\mathcal{ML}} = \mathcal{E}_{\lambda_2}^{\mathcal{ML}}$. However, the converse is not true in general.

 \begin{lemma}\label{lma3.8}
 Let $\lambda$ be a minimal measured lamination and $\mu$ be any measured lamination. Then $\mathcal{E}_{\lambda}^{\mathcal{ML}} = \mathcal{E}_{\mu}^{\mathcal{ML}}$ if and only if $\mu = \lambda' + c_1\gamma_1+ c_2 \gamma_2 + \dots + c_n \gamma_n$; where $\lambda'$ is a measured lamination supported on $|\lambda|$, $\gamma_{i}$ is a boundary curve of the subsurface  filled by $\lambda$ and $c_{i} \geq 0$ for all $i$.
 
 Consequently, if $\lambda$ and  $\mu$ are minimal measured laminations, then $\mathcal{E}_{\lambda}^{\mathcal{ML}} = \mathcal{E}_{\mu}^{\mathcal{ML}}$ if and only if $|\lambda| = |\mu|$.
 \end{lemma}

 \begin{proof} If $\lambda$ is a measure supported on a simple closed curve, then by \cref{lma3.6}
 
 $$\mathcal{E}_{\lambda}^{\mathcal{ML}}  = \{c\lambda : c \in \mathbb{R}_{+}\}= \mathcal{E}_{\mu}^{\mathcal{ML}}$$ if and only if $|\lambda| = |\mu|.$

 Now, assume $\lambda$ is a minimal lamination that fills a subsurface $S'$ of $S$ and $\mu = \lambda' + c_1\gamma_1+ c_2 \gamma_2 + \dots + c_n \gamma_n$; where $\gamma_{i}'s$ are the boundary curves of $S'$ and $|\lambda'| = |\lambda|$. 
 For any measured lamination $\delta$,
 $$
     i(\mu, \delta) = i(\lambda', \delta) + c_{1}i(\gamma_{1}, \delta) + \dots + c_{n}i(\gamma_n, \delta).
$$
 Any measured lamination that has zero intersection with $\mu$ must have a zero intersection with $\lambda'$, and therefore a zero intersection with $\lambda$. If $\delta$ belonged to $\mathcal{E}_{\mu}^{\mathcal{ML}}$, then for any measured lamination $c$,
 $$i(\delta, c) \neq 0 \Rightarrow i(\mu, c) \neq 0.$$
 This would mean at least one of the summand in 
 $$i(\mu, c) = i(\lambda', c) + c_{1}i(\gamma_{1}, c) + \dots + c_{n}i(\gamma_n, c)$$
 
 must be non zero. Since  $\lambda'$ and  the peripheral curves $\gamma_{i}$ all belong in $\mathcal{E}_{\lambda}^{\mathcal{ML}}$, any one of those summand being non zero implies $i( \lambda, c) \neq 0$. This gives us,
 
 $$\mathcal{E}_{\mu}^{\mathcal{ML}}  \subseteq \mathcal{E}_{\lambda}^{\mathcal{ML}}$$

 To see the other inclusion, observe that $\mathcal{E}_{\mu}^{\mathcal{ML}}$ contains measured laminations that are supported on its minimal components, $|\lambda|$, $|\gamma_1|$, $|\gamma_{2}|$, $\dots$ and $|\gamma_{n}|$. But, from \cref{lma3.7} we get 
 $$\mathcal{E}_{\lambda}^{\mathcal{ML}} \subseteq \mathcal{E}_{\mu}^{\mathcal{ML}} .$$

 To see the converse, assume that $\mathcal{E}_{\mu}^{\mathcal{ML}}$ = $\mathcal{E}_{\lambda}^{\mathcal{ML}}.$ The lamination $\mu$ cannot have a minimal component $\lambda_1$, that fills some subsurface of $S$ but is not supported in $|\lambda|$. If it does, then $\lambda_1$ will be contained in $\mathcal{E}_{\mu}^{\mathcal{ML}}$, but not in $\mathcal{E}_{\lambda}^{\mathcal{ML}}$. For the same reason, $|\mu|$ cannot contain  any simple closed curves that are  not a boundary curve of the subsurface $S'$. This means that $|\lambda|$ and $|\mu|$ can only differ by the simple closed curves that form the boundary curves of the subsurface filled by $\lambda$.
\end{proof}

\begin{remark}
 The construction $\mathcal{E}_{\lambda}$ can be generalized even further to the setting of geodesic currents. We will denote it by $\mathcal{E}^{c}_{\lambda}$, and it is defined as follows: If $\lambda$ is a geodesic current 
\[
\mathcal{E}^{c}_{\lambda} := 
\begin{cases} 
      \mu \in \mathscr{C} : & i(\lambda, \mu) = 0  \text{ and } \\
      
      & i(\lambda, k) \neq 0 \text{ for every geodesic current } k \text{ with } i(\mu, k) \neq 0.
\end{cases}
\]
 It can be shown that $\mathcal{E}^{c}_{\lambda}$ only consists of measured laminations. In fact, for a minimal measured lamination $\lambda$, the sets $\mathcal{E}^{c}_{\lambda}$ and $\mathcal{E}^{\mathcal{ML}}_{\lambda}$ are the same. For a closed curve $\gamma$, the sets $\mathcal{E}^{c}_{\gamma}$ and $\mathcal{E}^{\mathcal{ML}}_{\gamma}$ are the same as well.
For this reason, replacing $\mathcal{E}_{\lambda}^{\mathcal{ML}}$ by $\mathcal{E}^{c}_{\lambda}$ and $\mathcal{E}_{\mu}^{\mathcal{ML}}$ by $\mathcal{E}^{c}_{\mu}$ in the statement of \cref{lma3.8} still gives us a true statement. The proof of this statement follows closely to the proof of \cref{lma3.8}.
 \end{remark}

In the following lemma we prove that for any measured lamination $\lambda$ the operations $\mathcal{E}_{-}^{\mathcal{ML}}$ and $\phi$ commute for each $\phi \in Aut(\mathcal{ML})$

\begin{lemma} \label{lma3.10}
Let $\lambda$ be a measured lamination and $\phi \in Aut(\mathcal{ML})$. Then $\phi(\mathcal{E}_{\lambda}^{\mathcal{ML}}) = \mathcal{E}_{\phi(\lambda)}^{\mathcal{ML}}$.
\end{lemma}

\begin{proof} Let $\mu$ be a measured lamination. It is enough to prove that $\mu \in \mathcal{E}_{\phi(\lambda)}^{\mathcal{ML}} $ if and only if $  \phi^{-1}(\mu) \in \mathcal{E}_{\lambda}^{\mathcal{ML}}$.

 Observe that, $i(\mu, \phi(\lambda)) = 0 \Leftrightarrow  i(\phi^{-1}(\mu), \lambda) = 0.$
 
Let $\mu \in \mathcal{E}_{\phi(\lambda)}^{\mathcal{ML}}$. This gives, $i(\mu, c') \neq 0 \Rightarrow i(\phi(\lambda), c') \neq 0$ for any measured lamiantion $c'.$ But now,
\begin{alignat*}{2}
   & i(\phi^{-1}(\mu), c') && \neq 0\\ 
  \Rightarrow  \hspace{0.2cm}  &   i(\mu, \phi(c')) && \neq 0\\
  \Rightarrow  \hspace{0.2cm}  &    i(\phi(\lambda), \phi(c')) && \neq 0\\
  \Rightarrow  \hspace{0.2cm}  &    i(\lambda, c') && \neq 0   
\end{alignat*}
So, $\phi^{-1}(\mu) \in \mathcal{E}_{\lambda}^{\mathcal{ML}}$. A very similar argument will give the other implication, that is, if $i(\phi^{-1}(\mu), c') \neq 0 \Rightarrow i(\lambda, c') \neq 0$, then $i(\mu, c') \neq 0 \Rightarrow i(\phi(\lambda), c') \neq 0.$
This gives us, $\mu \in \mathcal{E}_{\phi(\lambda)}^{\mathcal{ML}} \Leftrightarrow \phi^{-1}(\mu) \in \mathcal{E}_{\lambda}^{\mathcal{ML}} \Leftrightarrow \mu \in \phi(\mathcal{E}_{\lambda}^{\mathcal{ML}}) .$
\end{proof}

\begin{remark}
    The same proof can be used to show that for a geodesic current  $\lambda$ and $\phi \in Aut(\mathscr{C})$, $\phi(\mathcal{E}^{c}_{\lambda}) = \mathcal{E}^{c}_{\phi(\lambda)}$.
\end{remark}

\begin{lemma}\label{lma3.12}
Automorphisms of measured laminations that preserve the intersection form map minimal measured laminations to minimal measured laminations.
\end{lemma}

\begin{proof} 
Let $\lambda$ be a minimal measured lamination on $S$ and let  $\phi \in Aut(\mathcal{ML})$.
Let 
$$\phi(\lambda) = \lambda_1 + \lambda_2 + \dots + \lambda_n$$
be the decomposition of the measured lamination $\phi(\lambda)$ into minimal measured laminations. Since $\phi^{-1} \in Aut(\mathcal{ML})$ is linear, we get
\[
\lambda = \phi^{-1}(\lambda_1) + \phi^{-1}(\lambda_2) + \dots + \phi^{-1}(\lambda_n)  
\]
Thus, support $|\phi^{-1}(\lambda_i)|$ is contained in $|\lambda|$ for all $i$. But since $\lambda$ is minimal $|\phi^{-1}(\lambda_i)| = |\lambda|$ for all $i$. Therefore, for any two minimal components $\lambda_{i}$ and $\lambda_{j}$
\begin{alignat*}{2}
    & \hspace{0.3cm}  |\phi^{-1}(\lambda_{i})| &&= |\phi^{-1}(\lambda_{j})|\\
    \Rightarrow  \hspace{0.2cm}  &\hspace{0.6cm}  \mathcal{E}_{\phi^{-1}(\lambda_{i})}^{\mathcal{ML}} &&= \mathcal{E}_{\phi^{-1}(\lambda_{j})}^{\mathcal{ML}}\\
    \Rightarrow  \hspace{0.2cm}  & \phi( \mathcal{E}_{\phi^{-1}(\lambda_{i})}^{\mathcal{ML}}) &&= \phi(\mathcal{E}_{\phi^{-1}(\lambda_{j})}^{\mathcal{ML}})\\
    \Rightarrow  \hspace{0.2cm}  &\hspace{1.4cm}   \mathcal{E}_{\lambda_{i}}^{\mathcal{ML}} &&= \mathcal{E}_{\lambda_{j}}^{\mathcal{ML}}\\
\end{alignat*}
Since $\lambda_{i}$ and $\lambda_{j}$ are minimal laminations, by ~\cref{lma3.8} we can conclude that $|\lambda_{i}| = |\lambda_{j}|$. Therefore, all the minimal components of $\phi(\lambda)$ have the same support and hence $\phi(\lambda)$ is a minimal measured lamination. 
\end{proof}

\begin{remark}
     Any $\phi \in Aut(\mathscr{C})$ also maps minimal measured laminations to minimal measured laminations. To see this replace the set $\mathcal{E}_{-}^{\mathcal{ML}}$ in the above proof with $\mathcal{E}^{c}_{-}$.

    A similar alteration to the proof of~\cref{prop3.2} will show that any $\phi \in Aut(\mathscr{C})$ maps simple closed curves to weighted simple closed curves. 
\end{remark}
\textbf{Proof of~\cref{prop3.2}:}

\begin{proof}
Let $\gamma$ be a simple closed curve on $S$. Then, 
\[
\mathcal{E}_{\gamma}^{\mathcal{ML}} = \{c\gamma: c \in \mathbb{R}_{+}\}.
\]
This implies, 
\[
    \phi(\mathcal{E}_{\gamma}^{\mathcal{ML}}) = \{c\phi(\gamma): c \in \mathbb{R}_{+}\}
\]
and so,
\begin{equation}
    \mathcal{E}_{\phi(\gamma)}^{\mathcal{ML}} = \{c\phi(\gamma): c \in \mathbb{R}_{+}\}.
\end{equation}

From \cref{lma3.12}, we know that $\phi(\gamma)$ is minimal. By \cref{lma3.6} and \cref{lma3.7}, $\mathcal{E}_{\phi(\gamma)}^{\mathcal{ML}} = \{c\phi(\gamma): c \in \mathbb{R}_{+}\}$ if and only if $|\phi(\gamma)|$ satisfies one of the following two scenarios:
\begin{enumerate}
    \item $|\phi(\gamma)|$ is a simple closed curve.
    \item $|\phi(\gamma)|$ is a uniquely ergodic minimal measured lamination that fills the entire surface $S$.
\end{enumerate}
Assume case $2$. Now, consider any simple closed curve $\alpha$ in $S$ such that $i(\gamma, \alpha) = 0$. This gives us $i(\phi(\gamma), \phi(\alpha)) = 0$. Since $\phi(\lambda)$ is a minimal lamination that fills all of $S$, \cref{lma3.5} implies that $|\phi(\gamma)| = |\phi(\alpha)|$.

Let us now pick another simple closed curve $\beta$, such that $i(\gamma, \beta) = 0$ and $i(\alpha, \beta) \neq 0$. This implies, $i(\phi(\gamma), \phi(\beta)) = 0$ and $i(\phi(\alpha), \phi(\beta)) \neq 0$. But earlier we concluded that $|\phi(\gamma)| = |\phi(\alpha)|$. By  ~\cref{lma3.4}, this is a contradiction.

The only possible case is $1$, $|\phi(\gamma)|$ is a simple closed curve. Therefore, $\phi(\gamma)$ is a weighted simple closed curve.
\end{proof}

\subsection{Preserving weights of simple closed curves}

We can now go one step further, and show that any $\phi \in Aut(\mathcal{ML})$ also preserves the weights of simple closed curves. That is, if  $\phi \in Aut(\mathcal{ML})$, then $\phi$ maps a simple closed curve with unit weight to a simple closed curve with unit weight. In order to prove this, we consider the action of $\phi$ on the curve complex on $S$.

\begin{lemma}\label{lemma3.9}
For any $\phi \in Aut(\mathcal{ML})$, we can find $\phi^{\prime} \in Aut(\mathcal{ML})$ so that for any simple closed curve $\gamma$, there is some $k$ so that $\phi^{\prime} \circ \phi(\gamma) = k \gamma$.
\end{lemma}

\begin{proof} 
Let $\left[ \gamma \right]$ denote the vertex in the curve complex that corresponds to the simple closed curve $\gamma$. 
Consider a $\phi \in Aut(\mathcal{ML})$. Let's say $\phi(\gamma) = k \gamma'$, where $k$ depends on the choice of $\gamma$. We define the action of $\phi$ on the curve complex by
\[
\phi^{*}(\left[ \gamma \right]) := \left[ \gamma' \right].
\]
This action is well defined as $\phi$ preserves disjointedness. If $i(\gamma_1, \gamma_2) = 0$, then $i(\phi(\gamma_1), \phi(\gamma_2)) = 0$.

By Ivanov's theorem we can find an $f \in Mod^{\pm}(S)$ such that its action on the curve complex agrees with $\phi^*$. 
That is
\[
f^*([\gamma]) = \phi^*([\gamma]) \text{ , for any closed curve } \gamma.
\]

But then $(f^{-1})^{*}(\left[ \gamma' \right]) = \left[ \gamma \right]$. This means that the action of $f^{-1}$ on $\mathcal{ML}(S)$ maps $\gamma'$ to $\gamma$. Let us use $\phi'$ to denote the action of $f^{-1}$ on $\mathcal{ML}(S)$. We have,
$$\phi'(\gamma') = \gamma$$
The map $\phi'$ satisfies the desired relation 
$$\phi' \circ \phi (\gamma) = \phi' (k \gamma') = k\gamma$$ \end{proof}

\begin{proposition}\label{prop3.17} Automorphisms of measured laminations that preserve the intersection form map simple closed curves to simple closed curves.
\end{proposition}
\begin{proof}
Let $\phi \in Aut(\mathcal{ML})$. From~\cref{prop3.2}, we know that $\phi$ maps simple closed curves to weighted simple simple closed curves. By \cref{lemma3.9}, we can find a $\phi' \in Aut(\mathcal{ML})$ such that $\phi' \circ \phi(\gamma) = k \gamma$ for all simple closed curves $\gamma$, and $k$ depending on $\gamma$. We will now show that the weight k has to be one for any simple closed curve.

Fix a simple closed curve $\gamma_{1}$ on our surface. Find simple closed curves $\gamma_2,\gamma_3$ so that $\gamma_1, \gamma_2$ and $\gamma_3$ all pairwise intersect. All of these can be summarised by the following notations

\[
\phi' \circ \phi(\gamma_{i}) = k_{i} \gamma_{i} \text{ for all }i \in \{ 1,2,3 \}
\]
and 
\[
i(\gamma_{i}, \gamma_{j}) \neq 0 \text{ for all }i,j \in \{1,2,3\}
\]
Observe that 
\[
    i(\gamma_{i}, \gamma_{j}) = i(\phi' \circ \phi (\gamma_{i}),  \phi' \circ \phi(\gamma_{j})) = i(k_{i} \gamma_{i}, k_{j} \gamma_{j}) =  k_{i} k_{j} \cdot i (\gamma_{i}, \gamma_{j})
\]

for all $i, j = 1,2,3$. This implies 

\begin{alignat*}{1}
   &k_{1}k_{2} = k_{1}k_{3} = 1\\
   \Rightarrow \hspace{0.2cm}   & k_{1}(k_{2} - k_{3}) = 0\\
   \Rightarrow \hspace{0.2cm}   & k_{1} = 0 \text{ or } k_{2} = k_{3}\\
\end{alignat*}
Since $k_{1}$ divides one, $k_{1}$ does not equal $0$. This gives us $k_{2} = k_{3}$. Furthermore, $k_{2}k_{3} = 1$ and the weights $k_{2}$ and $k_{3}$ are both positive. This proves that $k_{2} = k_{3} = 1$ and therefore $k_{1} = 1$.

Since the choice of our simple closed curve $\gamma_{1}$ was arbitrary, it follows that $\phi' \circ \phi(\gamma) = \gamma$ for all simple closed curve $\gamma$. The proof of \cref{lemma3.9} establishes that $\phi'$ is induced by an element belonging to $Mod^{\pm}(S)$. Consequently, this implies that $\phi(\gamma) = \gamma$ for all simple closed curve $\gamma$.
\end{proof}

\begin{remark} \label{rem3.18} If we replace $\phi \in Aut(\mathcal{ML})$ with $\phi \in Aut(\mathscr{C})$ in the proofs of \cref{lemma3.9} and \cref{prop3.17} we will see that the statements hold for automorphisms of currents as well. 
Any $\phi \in Aut(\mathscr{C})$ maps simple closed curves to simple closed curves. 
\end{remark}

\section{Proof of the theorem 1.1}

\begin{reptheorem}{thm0.1}
Let $S_{g}$ be a closed, orientable, finite type surface of genus $g \geq 2$ and let $Aut(\mathcal{ML})$ denote the group of homeomorphisms on $\mathcal{ML}(S_{g})$ that preserve the intersection form. Then
\[
Aut(\mathcal{ML}) \cong Mod^{\displaystyle \pm} (S_{g})
\]
for all $g \neq 2$. For the surface of genus 2, 
\[
Aut(\mathcal{ML}) \cong \displaystyle{Mod^{ \pm} (S_{2})/H}
\]
Where, $H$ is the order two subgroup generated by the hyperelliptic involution.   
\end{reptheorem}

\begin{proof} 
Let us denote the automorphism group of the curve complex on a surface $S_{g}$ by $Aut(CC)$. From Ivanov's theorem we know that $Aut(CC) \cong Mod^{\displaystyle \pm} (S_{g})$, when $g \neq 2$. And for the surface of genus 2, $Aut(CC) \cong \displaystyle{Mod^{ \pm} (S_{2})/H}$. We will show that 
$Aut(\mathcal{ML}) \cong Aut(CC)$.

Consider the map $\psi: Aut(\mathcal{ML}) \rightarrow Aut(CC)$ that maps automorphisms on $\mathcal{ML}$ to its action on the curve complex.
\[
\psi(\phi) = \phi^{*}
\]

To see $\psi$ is injective, consider any element $\phi$  in the kernel of $\psi$. By \cref{prop3.17} any such $\phi$ fixes all the simple closed curves. Because $\phi$ is linear by \cref{prop3.1}, it also fixes all the weighted multicurves. But, weighted multicurves are dense in $Aut(\mathcal{ML})$. Therefore, $\phi$ is identity on $\mathcal{ML}$, and consequently $\psi$ is injective.

Observe that, $Mod^{\pm}(S_{g})$ has a natural action on $\mathcal{ML}(S_{g})$ that also preserves the intersections form. For every $f \in Mod^{\pm}(S_{g})$, we can find a $\widetilde{f} \in Aut(\mathcal{ML})$. This shows that $\psi$ is surjective.
\end{proof}

\section{Obstruction when generalizing to currents}\label{Sec5}
For a surface $S$ of genus at least 3, \cref{thm0.1} gives us a surjection from $Aut(\mathscr{C})$ to $Mod^{\pm}(S)$.
\[ 
Aut(\mathscr{C}) \xtwoheadrightarrow{f} Aut(\mathcal{ML}) \cong Aut(CC) \cong Mod^{\pm}(S)
\]
We are interested to know if the map $f$ is injective. The kernel of $f$ consists of maps in $Aut(\mathscr{C})$ that induces an identity map on $Aut(CC)$. In other words, any $\phi$ in the kernel of $f$ fixes all simple closed curves.
Now, if $f$ is injective, then the kernel of $f$ will only contain the identity map on currents. This means that $f$ being injective will imply the following statement: Any automorphism of currents that preserves the intersection form and fixes all simple closed curves is the identity on currents. 

 Clearly, if $\phi$ is in the kernel of $f$ then it preserves the simple marked length spectrum of any current $\mu$. However, this is not enough to prove $f$ is injective. In fact, in this section we construct an infinite family of closed curves that have the same self intersection number and simple marked length spectrum. Let $\gamma$ and $\delta$ be two such closed curves. It is not immediate why there cannot exist a $\phi \in ker(f)$ that maps $\gamma$ to $\delta$. The remaining part of this section will be focused on constructing these examples. We will also remark about some other properties of these curves.

\begin{reptheorem}{thm0.3} 
Let $S$ be a surface of genus at least 2.
We can find infinitely many pairs of closed curves $\gamma_n$ and $\gamma_n'$ on $S$ such that $\gamma_n$ and $\gamma_n'$ have the same self intersection number and the same simple marked length spectra. 
\end{reptheorem}

\begin{proof}

Let $S$ be closed surface of genus at least 2. Let $P$ be a pair of pants on $S$ and let $\partial P$ be its boundary consisting of oriented curves $x, y$ and $z$ labeled below. We will use $\overline{x}$, $\overline{y}$, $\overline{z}$ to denote the curves that go along the boundary curves $x$, $y$ and $z$ but have the opposite orientation. We consider on $P$ the family of closed curves as given in Figure 3.

For each triple $(a,b,c)$ where $a$, $b$, $c$ are positive integers with $a \geq b \geq c$ we consider the closed curve $\gamma_{(a,b,c)}$. Let $x_0$ be the point on the curve at labelled in Figure 3. The path of the curve $\gamma_{(a,b,c)}$ is described as follows: $(1)$ Starting at $x_0$ it takes one full twist along $\overline{z}$, $(2)$ $c$ half twists about $y$, $(3)$ one full twist along $\overline{x}$, $(4)$ $b$ half twists about $y$, $(5)$ one full twist along $\overline{x}$, $(6)$ $a$ half twists about $y$, and then the curve closes. Any two curves in this family are identical, except for the number of half twists $a$, $b$ and $c$. Note that the integer $c$ is always odd since the curve follows different boundary components before and after making the $c$ half twists along $y$. For a similar reason $a$ is also odd. But, $b$ is even as the curve follows along the same curve $\overline{x}$ before and after making the $b$ half twists about $y$.

\begin{figure}[h]
    \centering
    \includegraphics[width=\linewidth]{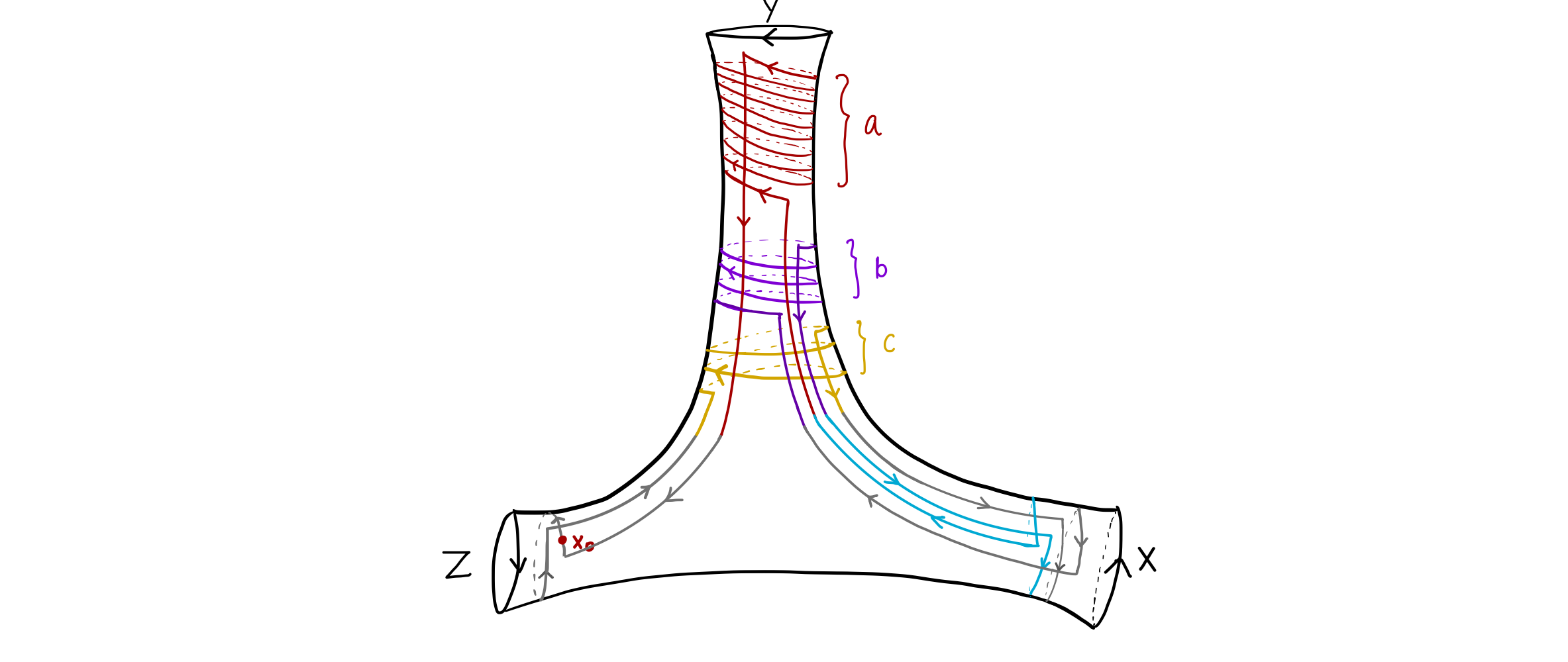}
    \caption{Closed curve with number of half twists $a =19$, $b = 8$ and $c = 5$. The figure shows the curve $\gamma_{(19,8,5)}$ intersecting itself minimally.}
\end{figure}

The curve $\gamma_{(a,b,c)}$ shown in figure 3 has no bigons. As a result the curve is in the minimal position that realizes its self intersection as long as $a \geq b \geq c$ \cite{Primer}. If we count the number of self intersections we get:
\begin{equation}\label{eq3}
  i(\gamma_{(a,b,c)}, \gamma_{(a,b,c)}) = \displaystyle{ \Big( \frac{a-1}{2} \Big) + 3 \Big(\frac{b-2}{2} \Big) + 5 \Big( \frac{c-1}{2} \Big) + 6 }  
\end{equation}

The curves $\gamma_{(a,b,c)}$ and $\gamma_{(a',b',c')}$ that have the same self intersection number satisfy:

$$\displaystyle{ \Big( \frac{a-1}{2} \Big) + 3 \Big(\frac{b-2}{2} \Big) + 5 \Big( \frac{c-1}{2} \Big) + 6 } = \displaystyle{ \Big( \frac{a'-1}{2} \Big) + 3 \Big(\frac{b'-2}{2} \Big) + 5 \Big( \frac{c'-1}{2} \Big) + 6 }$$

This is equivalent to 
\begin{equation}\label{sumeq}
    (a - a') + 3(b - b') + 5(c - c') = 0.
\end{equation}

In what follows, we will show that for $\gamma_{(a,b,c)}$ and $\gamma_{(a',b',c')}$ to have the same simple marked length spectra it is enough that the following equation is satisfied:
\begin{equation}\label{eqsa}
     (a - a') + (b - b') + (c - c') = 0
\end{equation}

A simple closed curve in $S$ intersects $\gamma_{(a,b,c)}$ or $\gamma_{(a',b',c')}$ if and only if it passes through $P$. Any such simple closed curve will intersect $P$ as a simple essential arc with its end points in $\partial P$. For any essential arc $s$ in $P$ with its end points in $\partial P$, we want \cite[Proposition 62]{CL}:
\begin{equation*}
i(s, \gamma_{(a,b,c)}) = i(s,\gamma_{(a',b',c')})  
\end{equation*}

Here, we let $i(\cdot, \cdot)$ denote the minimum intersection between $\gamma_{(a,b,c)}$ and any arc in the free homotopy class of $s$, where the endpoints of $s$ to remain in $\partial P$ throughout the homotopy; the end points of $s$ need not be fixed. Up to homotopy we have six such essential simple arcs on $P$. Since the curves $\gamma_{(a,b,c)}$ and $\gamma_{(a',b',c')}$ are identical except for the number of half twists about $Y$, the three essential arcs with end points only in $x$ and $z$ will intersect $\gamma_{(a,b,c)}$ and $\gamma_{(a',b',c')}$ the same number of times.

\begin{figure}[h]
    \centering
    \includegraphics[width=\linewidth]{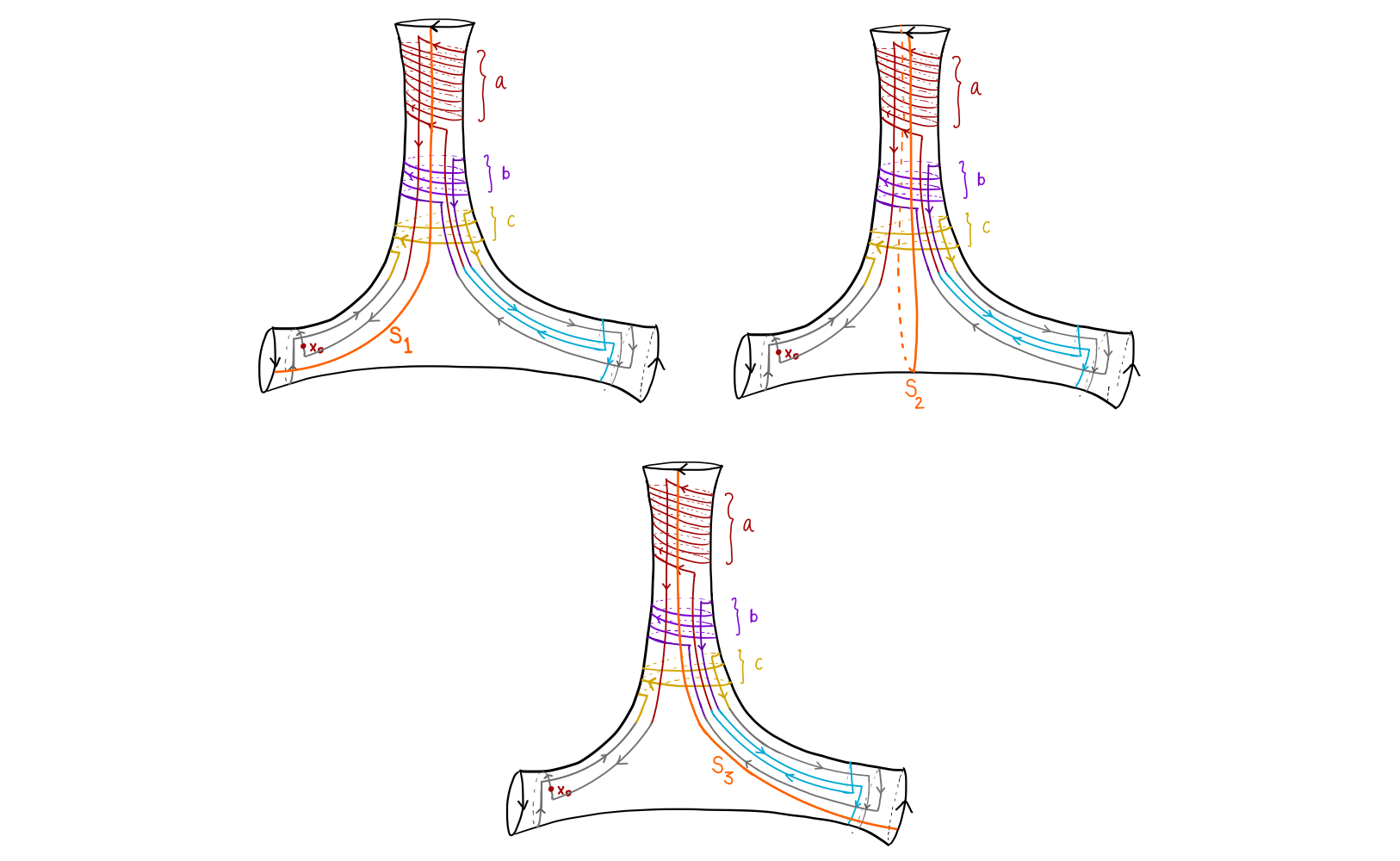}
    \caption{Arcs $s_1$, $s_2$ and $s_3$ are homotoped to intersect $\gamma_{(a,b,c)}$ minimal number of times.}
    \label{fig:enter-label}
\end{figure}

Figure 4 shows the other three essential simple arcs on $P$ up to homotopy. It is straightforward to see that all the three arcs $s_1$, $s_2$ and $s_3$ intersect $\gamma_{a,b,c}$ minimally in Figure 4.

Counting the number intersections we get 
\begin{equation}
    \displaystyle{i(\gamma_{(a,b,c)}, s_1) = \Big( \frac{a-1}{2} \Big) + \Big( \frac{b-2}{2} \Big) + \Big( \frac{c-1}{2} \Big) + 3 }
\end{equation}
\begin{equation}
    \displaystyle{i(\gamma_{(a,b,c)}, s_2) = a + b + c} 
\end{equation}
\begin{equation}
   \displaystyle{ i(\gamma_{(a,b,c)}, s_3) = \Big( \frac{a-1}{2} \Big) + \Big( \frac{b-2}{2} \Big) + \Big( \frac{c-1}{2} \Big) + 4 }
\end{equation}

Replacing $a$, $b$, $c$ in the above equations by $a'$, $b'$ and $c'$ respectively give us the intersection number between $\gamma_{(a,',b',c')}$ and the arcs. The condition 
\[
i(s_j, \gamma_{(a,b,c)}) = i(s_j,\gamma_{(a',b',c')})  
\]
for $j = 1, 2$ and $3$ is equivalent to \cref{eqsa} . 

For $k$ an even positive integer and let $t$ an odd integer greater than or equal to 3, we set
   \begin{align}\label{eq14}
       &c = t &c'&= k + t\\
       &b = 4k +t+1  &b'&= 2k + t+ 1 \nonumber \\
       &a = 7k +t  &a'&= 8k + t \nonumber
   \end{align}
   Then $(a,b,c)$ and $(a',b',c')$ satisfy equations (3) and (4). They also satisfy $a \geq b \geq c$ and $a' \geq b' \geq c'$. The conditions on $k$ and $t$ make certain that $a$, $a'$, $c$, $c'$ are odd and $b$, $b'$ are even. This gives us a infinite collection of curves $\gamma_{(a,b,c)}$ and $\gamma_{(a',b',c')}$ have the same self intersection number and simple marked length spectra. 
\end{proof}
\begin{remark}
    The solution set that satisfies \cref{sumeq} and \cref{eqsa} is 4-dimensional. The collection of curves $\gamma_{(a,b,c)}$ and $\gamma_{(a',b',c')}$ that have the same self intersection number and simple marked length spectra is bigger than the one described in \cref{eq14}. The two parameter family described in \cref{eq14} makes computing examples easier. 
\end{remark}

\begin{remark}The curves $\gamma_{(a,b,c)}$ and $\gamma_{(a',b',c')}$ are not hyperbolically equivalent. To see this consider the fundamental group of $P$. This is the free group $F(x,y)$ where $x$ and $y$ represent the loops that go around the boundary components $X$ and $Y$ respectively. For any $a, b, c$, the word in the fundamental group representing the curve $\gamma_{(a,b,c)}$ is:
    \begin{equation}\label{eqpi1}
        \displaystyle{yxy^{\frac{c+1}{2}}x^{-1}y^{\frac{b}{2}}x^{-1}y^{\frac{a-1}{2}}} 
    \end{equation}
For any $k$ and $t$ we will always have $a \neq a'$, $b \neq b'$ and $c \neq c'$. The words in the fundamental group representing the curves $\gamma_{(a,b,c)}$  and $\gamma_{(a',b',c')}$ will have different exponents of $y$.  By Theorem 1.2 in \cite{Hor} the two curves in the same pair are not trace equivalent. By the work on Leininger the curves cannot be hyperbolically equivalent on $S$ \cite{CL}. 
\end{remark}

\begin{remark} These examples on a pair of pants allow us to define quite complicated pairs of examples on arbitrary surfaces. Let $S'$ be a surface of genus at least two. One can define $\pi_1$-injective, immersion maps $f: P \rightarrow S'$ that takes boundary curves $X$ and $Y$ of $P$ to distinct closed curves in $S'$ as given in figure 5. Any simple closed curve $\alpha$ in $S'$ intersects $f(P)$ as a simple arc. The pre-image $f^{-1}(\alpha)$ is a disjoint union of simple arcs. These arcs intersect $\gamma_{(a,b,c)}$ and $\gamma_{(a',b',c')}$ the same number of times. The curves $f(\gamma_{(a,b,c)})$ and $f(\gamma_{(a',b',c')})$ on $S'$ share the same self intersection number. 

It is known that length equivalence is a property that is preserved under continuous map \cite{CL}. The curves $f(\gamma_{(a,b,c)})$ and $f(\gamma_{(a',b',c')})$ are not length equivalent. This allows us to construct pairs of filling closed curves on higher genus surfaces that not length equivalent but are simple intersection equivalent. 
\end{remark}

\begin{figure}[h]
    \centering
    \includegraphics[width=\linewidth]{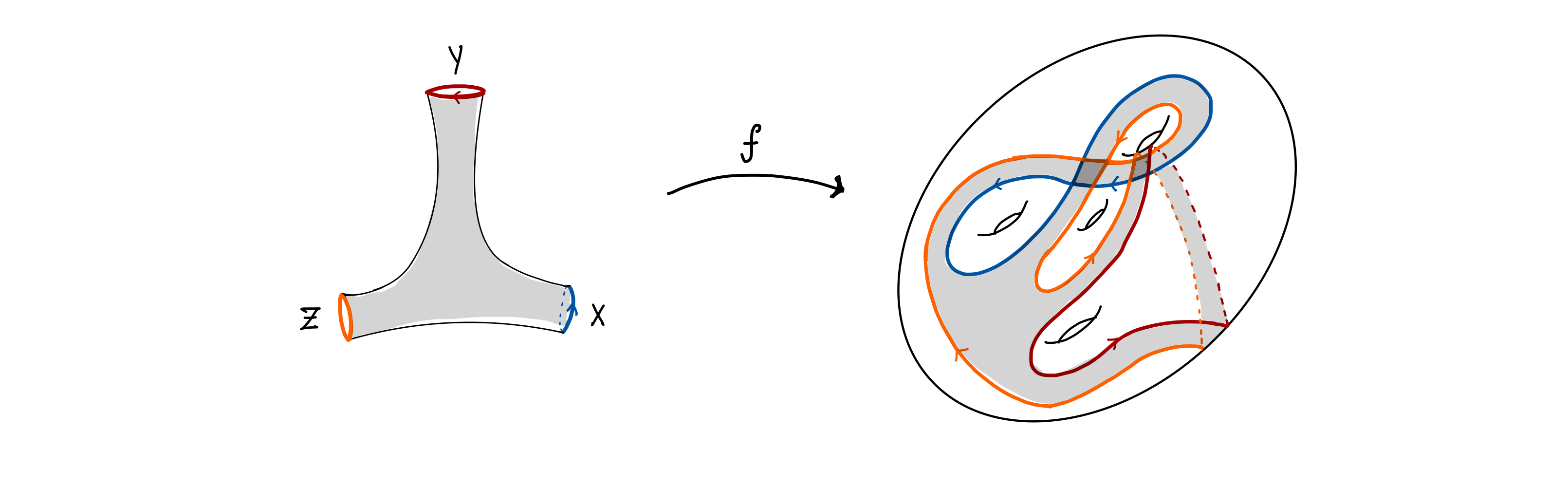}
    \includegraphics[width=\linewidth]{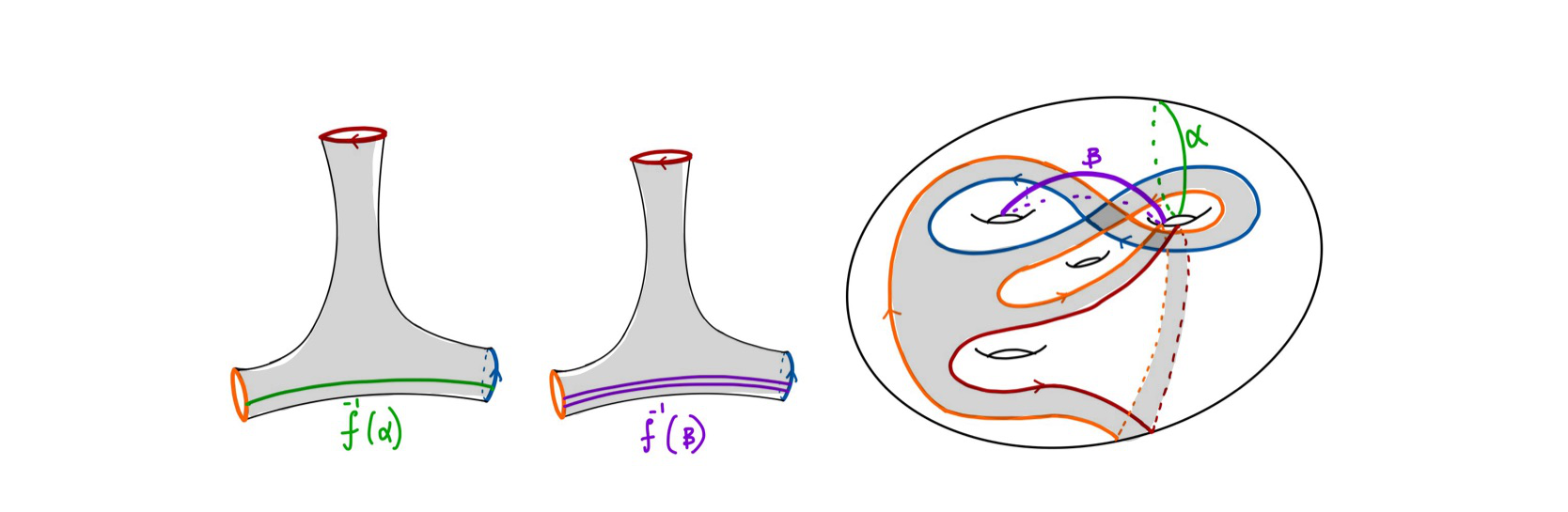}
    \caption{(Top) Shaded region is the image of the immersion $f: P \rightarrow S'$. The boundary curves $X$, $Y$, $Z$ are mapped to closed curves on $S'$ of the respective color. (Bottom right) $\alpha$ and $\beta$ are simple closed curves on $S'$ that intersect $f(P)$. (Bottom left) The arcs in $P$ are inverse images of $\alpha$ and $\beta$ in $P$.}
    \label{fig:enter-label}
\end{figure}

\section{Related questions}\label{Sec6}

 From the remarks in section 3, we know that some results for $\mathcal{ML}$ also holds true in the general setting of currents. In particular, we know that any $\phi \in Aut(\mathscr{C})$ is linear and takes simple closed curves to simple closed curves. This gives an action of $Aut(\mathscr{C})$ on the curve complex. However, we do not know whether this action is faithful. Because of this, it is not straightforward to see if $Aut(\mathscr{C})$ is $Mod^{\pm}(S)$.

In fact, it is not yet known whether maps in $Aut(\mathscr{C})$ take closed curves to closed curves. 

\begin{question}
    Do automorphisms of geodesic currents that preserve the intersection form map closed curves to closed curves? 
\end{question}

Ivanov's meta-conjecture has been verified for the automorphism group of different objects associated with a surface, such as the Teichmüller space, the arc complex, and \( k \)-graph \cite{Irm06, Di12, AACLOSX, CW2015, IK07}. Many of these proofs demonstrate that the automorphism group in question admits an action on the curve graph. Additionally, several of these objects including the Teichmüller space, arc complex, and the \( k \)-graph can be determined using simple closed curves. We know that closed curves cannot be determined by their intersections with simple closed curves. Moreover, simple closed curves are a very small and discrete subset of the space of currents. It would be interesting to understand the group \( \operatorname{Aut}(\mathscr{C}) \), and if this group is the mapping class group, it would be interesting to understand why maps of currents that preserve the intersection form can be determined by their action on simple closed curves.

\begin{question}
    For a closed orientable finite type surface $S$ of genus at least two, is $Aut(\mathscr{C}) \cong Mod^{\pm}(S)$?
\end{question}

\bibliographystyle{alpha} 
\bibliography{paperdraft}

\newcommand{\etalchar}[1]{$^{#1}$}
\begin{thebibliography}{AAC{\etalchar{+}}19}

\bibitem[AAC{\etalchar{+}}19]{AACLOSX}
Shuchi Agrawal, Tarik Aougab, Yassin Chandran, Marissa~Kawehi Loving, J.~Robert Oakley, Roberta Shapiro, and Yang Xiao.
\newblock Automorphisms of the k-curve graph.
\newblock {\em Michigan Mathematical Journal}, 2019.

\bibitem[AL17]{JACL}
Javier Aramayona and {Christopher J.} Leininger.
\newblock {\em Hyperbolic Structures on Surfaces and Geodesic Currents}, pages 111--149.
\newblock Advanced Courses in Mathematics - CRM Barcelona. Birkh{\"a}user, Switzerland, 2017.

\bibitem[And02]{Hor}
James Anderson.
\newblock Variations on a theme of horowitz.
\newblock 299, 02 2002.

\bibitem[BIPP17]{BIPP}
M.~Burger, Alessandra Iozzi, Anne Parreau, and Maria~Beatrice Pozzetti.
\newblock A structure theorem for geodesic currents and length spectrum compactifications.
\newblock {\em arXiv: Geometric Topology}, 2017.

\bibitem[BIPP21]{BIPP2}
M.~Burger, A.~Iozzi, A.~Parreau, and M.~B. Pozzetti.
\newblock Currents, systoles, and compactifications of character varieties.
\newblock {\em Proceedings of the London Mathematical Society}, 123(6):565--596, 2021.

\bibitem[BM19]{BM19}
Tara~E. Brendle and Dan Margalit.
\newblock Normal subgroups of mapping class groups and the metaconjecture of {I}vanov.
\newblock {\em J. Amer. Math. Soc.}, 32(4):1009--1070, 2019.

\bibitem[Bon86]{Bon86}
Francis Bonahon.
\newblock Bouts des vari\'{e}t\'{e}s hyperboliques de dimension {$3$}.
\newblock {\em Ann. of Math. (2)}, 124(1):71--158, 1986.

\bibitem[CB88]{CB88}
Andrew~J. Casson and Steven~A. Bleiler.
\newblock {\em Structure of geodesic laminations}, page 60–74.
\newblock London Mathematical Society Student Texts. Cambridge University Press, 1988.

\bibitem[Dis12]{Di12}
Valentina Disarlo.
\newblock Combinatorial rigidity of arc complexes.
\newblock {\em arXiv: Geometric Topology}, 2012.

\bibitem[DS03]{DS03}
Raquel D\'{\i}az and Caroline Series.
\newblock Limit points of lines of minima in {T}hurston's boundary of {T}eichm\"{u}ller space.
\newblock {\em Algebr. Geom. Topol.}, 3:207--234, 2003.

\bibitem[FM12]{Primer}
Benson Farb and Dan Margalit.
\newblock {\em A Primer on Mapping Class Groups (PMS-49)}.
\newblock Princeton University Press, 2012.

\bibitem[IK07]{IK07}
Elmas Irmak and Mustafa Korkmaz.
\newblock Automorphisms of the {H}atcher-{T}hurston complex.
\newblock {\em Israel J. Math.}, 162:183--196, 2007.

\bibitem[Irm06]{Irm06}
Elmas Irmak.
\newblock Complexes of nonseparating curves and mapping class groups.
\newblock {\em Michigan Math. J.}, 54(1):81--110, 2006.

\bibitem[Iva97]{Iva97}
Nikolai~V. Ivanov.
\newblock Automorphism of complexes of curves and of {T}eichm\"{u}ller spaces.
\newblock {\em Internat. Math. Res. Notices}, (14):651--666, 1997.

\bibitem[Iva06]{Iva06}
Nikolai~V. Ivanov.
\newblock Fifteen problems about the mapping class groups.
\newblock In {\em Problems on mapping class groups and related topics}, volume~74 of {\em Proc. Sympos. Pure Math.}, pages 71--80. Amer. Math. Soc., Providence, RI, 2006.

\bibitem[Kav19]{kavi2019}
Nithin Kavi.
\newblock Equivalence relations between closed curves on the pair of pants.
\newblock {\em arXiv preprint}, 2019.

\bibitem[Kon24]{Ker}
Hokuto Konno.
\newblock Dehn twists and the {N}ielsen realization problem for spin 4-manifolds.
\newblock {\em Algebr. Geom. Topol.}, 24(3):1739--1753, 2024.

\bibitem[Lei03]{CL}
Christopher~J. Leininger.
\newblock Equivalent curves in surfaces.
\newblock {\em Geometriae Dedicata}, 102:151--177, 2003.

\bibitem[Mar16]{Mar2016}
Bruno Martelli.
\newblock An introduction to geometric topology.
\newblock 2016.

\bibitem[OP18]{OHSHIKA2018899}
Ken'ichi Ohshika and Athanase Papadopoulos.
\newblock Homéomorphismes et nombre d'intersection.
\newblock {\em Comptes Rendus Mathematique}, 356(8):899--902, 2018.

\bibitem[Ota90]{Ot90}
Jean-Pierre Otal.
\newblock Le spectre marque des longueurs des surfaces a courbure negative.
\newblock {\em Annals of Mathematics}, 131(1):151--162, 1990.

\bibitem[PX23]{HpBx2023}
Hugo Parlier and Binbin Xu.
\newblock A topological viewpoint on curves via intersection.
\newblock {\em arXiv preprint}, 2023.

\bibitem[Wal15]{CW2015}
Cormac Walsh.
\newblock The horoboundary and isometry group of thurston's lipschitz metric.
\newblock {\em arXiv preprint}, 2015.

\end{thebibliography}

\end{document}